\documentclass[12pt,oneside,e-only]{amsart}

\usepackage{amssymb}
\usepackage{mathrsfs}
\usepackage{pstricks}
\usepackage{pst-node}
\usepackage{graphicx}

\date{}
\title{Word posets, complexity, and Coxeter groups}
\author{Matthew J. Samuel}
\address{Department of Mathematics, Rutgers University, 110
  Frelinghuysen Road, Piscataway, NJ 08854, USA}
\email{msamuel@math.rutgers.edu}

\begin{document}

\begin{abstract}
A monoid $M$ generated by a set $S$ of symbols can be described as the set of equivalence classes of finite words in $S$ under some relations that specify when some contiguous sequence of symbols can be replaced by another. If $a,b\in S$, a relation of the form $ab=ba$ is said to be a \emph{commutation relation}. Words that are equivalent using only a sequence of commutation relations are said to be in the same \emph{commutation class}. We introduce certain partially ordered sets that we call \emph{word posets} that capture the structure of commutation classes in monoids. The isomorphism classes of word posets are seen to be in bijective correspondence with the commutation classes, and we show that the linear extensions of the word poset correspond bijectively to the words in the commutation class, using which we demonstrate that enumerating the words in commutation classes of monoids is \#P-complete. We then apply the word posets to Coxeter groups. We show that the problem of enumerating the reduced words of elements of Coxeter groups is \#P-complete. We also demonstrate a method for finding the word posets for commutation classes of reduced words of an element, then use this to find a recursive formula for the number of commutation classes of reduced words.
\end{abstract}

\maketitle
\newtheorem{lemma}{Lemma}
\newtheorem{theorem}{Theorem}[section]
\newtheorem{proposition}[theorem]{Proposition}
\newtheorem{conjecture}[theorem]{Conjecture}
\newtheorem{corollary}{Corollary}
\theoremstyle{definition} \newtheorem*{defn}{}
\theoremstyle{definition} \newtheorem{example}{Example}

\newcommand{\boxrow}[3]{\SpecialCoor\psgrid[gridlabels=0pt,subgriddiv=1](#1,#2)(!#1 #3 add #2 -1 add)}

%\section*{Introduction}

\section{Word posets}

If $S$ is a set of symbols, a \emph{word} in $S$ is an element of the free monoid $S^{\ast}$, which is the set of all finite strings of symbols from $S$ with concatenation as the product and the empty string as the identity element. Any monoid $M$ generated by $S$ is a quotient of $S^{\ast}$ by some relations, which indicate when some contiguous sequence of symbols can be replaced by another, and $M$ is the set of these equivalence classes. A relation of the form $ab=ba$ is said to be a \emph{commutation relation}, and $a$ and $b$ are said to \emph{commute}. Two words equivalent using only a sequence of commutation relations are said to be in the same \emph{commutation class}.\\

\begin{example} \label{example:cclass}
Suppose we have a monoid generated by the symbols ${a}$, ${b}$, ${c}$, and ${d}$, satisfying the commutation relations
$${ab}={ba}\mathrm{,}$$
$$cd=dc\mathrm{,}$$
and
$${ad}={da}\mathrm{.}$$
Then the commutation class of ${abcd}$ consists of the words ${abcd}$, ${bacd}$, ${abdc}$, ${badc}$, and ${bdac}$.
\end{example}

Let $S$ be a set of generators, and let $\mathbf{w}\in S^{\ast}$ be a word in $S$. Denote by $\ell(\mathbf{w})$ the length of $\mathbf{w}$. We define a relation $\to$ on the set $[\ell(\mathbf{w})]=\{1,2,\ldots,\ell(\mathbf{w})\}$ by declaring that $i\to j$ if $i\leq j$ and either $\mathbf{w}_i$ and $\mathbf{w}_j$ are equal or do not commute. We then define a partial ordering $\leq_{\mathbf{w}}$ on $[\ell(\mathbf{w})]$ as the transitive closure of $\to$. Note that if $i\leq_{\mathbf{w}} j$, then $i\leq j$.\\

A poset $P$ with a function $s:P\to S$ will be called a \emph{word poset} if the following conditions are satisfied for all $u,v\in P$:\\
\begin{itemize}
\item[(a) ] If $s(u)$ and $s(v)$ are equal or do not commute, then $u$ and $v$ are comparable.\\

\item[(b) ] If $v$ covers $u$, then $s(u)$ and $s(v)$ either do not commute or are equal.\\
\end{itemize}

We can see that $\leq_{\mathbf{w}}$ together with $i\mapsto \mathbf{w}_i$ turns $[\ell(\mathbf{w})]$ into a word poset, and we will denote this poset by $P(\mathbf{w})$.

\begin{example}
For the word in Example \ref{example:cclass}, $P(\mathbf{w})$ can be illustrated by

$$\setlength{\arraycolsep}{1cm}
\psset{nodesep=3pt}
\psset{radius=3pt}
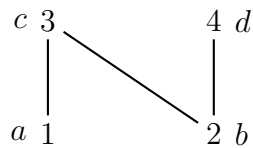
\begin{array}{cc}
\vspace{1cm}\rnode{C}{3}&\rnode{D}{4}\\
\rnode{A}{1}&\rnode{B}{2}
\end{array}
\ncline{A}{C}
\ncline{B}{C}
\ncline{B}{D}
\nput{180}{A}{{a}}
\nput{0}{B}{{b}}
\nput{180}{C}{{c}}
\nput{0}{D}{{d}}$$

\vspace{12pt}
With the same relations, if $\mathbf{w}={badbcd}$, then $P(\mathbf{w})$ is
$$\setlength{\arraycolsep}{1cm}
\psset{nodesep=3pt}
\psset{radius=3pt}
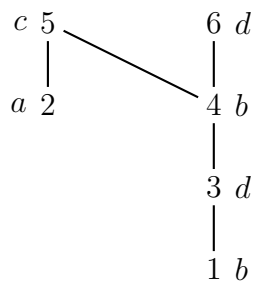
\begin{array}{cc}
\rnode{E1}{5}\vspace{0.6cm}&\rnode{F1}{6}\\
\vspace{0.6cm}\rnode{B1}{2}&\rnode{D1}{4}\\
\vspace{0.6cm}&\rnode{C1}{3}\\
\vspace{0.6cm}&\rnode{A1}{1}
\end{array}
\ncline{A1}{C1}
\ncline{B1}{E1}
\ncline{C1}{D1}
\ncline{D1}{E1}
\ncline{D1}{F1}
\nput{0}{A1}{{b}}
\nput{180}{B1}{{a}}
\nput{0}{C1}{{d}}
\nput{0}{D1}{{b}}
\nput{180}{E1}{{c}}
\nput{0}{F1}{{d}}
$$
\end{example}

\vspace{20pt}

An isomorphism of word posets $f:(P,s)\to (P',s')$ is a poset isomorphism such that $s'\circ f=s$.

\begin{lemma}
If $\mathbf{w}$ and $\mathbf{w}'$ are in the same commutation class, then there is an isomorphism $f:P(\mathbf{w})\to P(\mathbf{w}')$.
\end{lemma}
\begin{proof}
Suppose there is a $j$ such that $\mathbf{w}_i=\mathbf{w}'_i$ if $i\neq j,j+1$, $\mathbf{w}_j=\mathbf{w}'_{j+1}$, and $\mathbf{w}_{j+1}=\mathbf{w}'_j$. Then $\mathbf{w}_j$ and $\mathbf{w}_{j+1}$ commute, so an isomorphism is given by $f(i)=i$ if $i\neq j,j+1$, $f(j)=j+1$, and $f(j+1)=j$. If $\mathbf{w}$ and $\mathbf{w}'$ differ by more than one transposition, we may compose functions of the type just constructed to obtain our isomorphism $P(\mathbf{w})\to P(\mathbf{w}')$.
\end{proof}

\begin{theorem}
The words in the same commutation class as $\mathbf{w}$ are in bijective correspondence with the linear extensions of $P(\mathbf{w})$. If $e:P(\mathbf{w})\to[\ell(\mathbf{w})]$ is a linear extension, then the corresponding word $\mathbf{w}(e)$ is given by $\mathbf{w}(e)_i=\mathbf{w}_{e^{-1}(i)}$.
\end{theorem}
\begin{proof}
Write $k=\ell(\mathbf{w})$, and let $e$ be a linear extension of $P(\mathbf{w})$. Let $i$ be the index such that $e(i)=k$. If we define $f$ by $f(j)=e(j)$ if $j<i$, $f(j)=e(j+1)$ if $i\leq j<k$, and $f(k)=k$, then $f$ is also a linear extension of $P(\mathbf{w})$. By induction, $f$ restricted to $[k-1]$ corresponds to a word $\mathbf{w}''$ in the same commutation class as $\mathbf{w}'$ of length $k-1$ defined by $\mathbf{w}'_j=\mathbf{w}_j$ if $j<k$. Then $\mathbf{w}_k$ commutes with $\mathbf{w}'_j$ for all $j\geq i$, so $\hat{\mathbf{w}}$ defined by $\hat{\mathbf{w}}_j=\mathbf{w}''_j$ if $j<i$, $\hat{\mathbf{w}}_i=k$, and $\hat{\mathbf{w}}_j=\mathbf{w}''_{j-1}$ if $j>i$ is in the same commutation class as $\mathbf{w}$. But this is exactly $\mathbf{w}(f)$. Thus, $\mathbf{w}$ is a map from the set of linear extensions of $P(\mathbf{w})$ into the commutation class. If $\mathbf{w}'$ is in the same commutation class as $\mathbf{w}$, let $f:P(\mathbf{w})\to P(\mathbf{w}')$ be an isomorphism. Then $g$ is a linear extension of $P(\mathbf{w})$, and $\mathbf{w}'_{i}=\mathbf{w}_{g^{-1}(i)}$, so $\mathbf{w}'=\mathbf{w}(g)$. Thus, $\mathbf{w}$ is a surjection from the set of linear extensions onto the commutation class.

Next we need to see that $\mathbf{w}$ is injective. Suppose $e$ and $f$ are distinct linear extensions of $P(\mathbf{w})$, and let $k$ be the first index such that $e^{-1}(k)\neq f^{-1}(k)$. If $\mathbf{w}(e)=\mathbf{w}(f)$, then we must have that $\mathbf{w}_{e^{-1}(k)}=\mathbf{w}_{f^{-1}(k)}$. If $e^{-1}(k)<f^{-1}(k)$, then we must have that $e^{-1}(k)<_{\mathbf{w}} f^{-1}(k)$. But $f^{-1}(k)$ occurs earlier than $e^{-1}(k)$ in the linear extension $f$, which is a contradiction. Thus, $\mathbf{w}(e)\neq\mathbf{w}(f)$.
\end{proof}

\begin{theorem}
The correspondence $\mathbf{w}\mapsto P(\mathbf{w})$ establishes a bijection between commutation classes of words and isomorphism classes of word posets.
\end{theorem}
\begin{proof}
If $P$ is a word poset, let $e:P\to[|P|]$ be a linear extension, and define $\mathbf{w}_i=s(e^{-1}(i))$. Then $e$ is an isomorphism from $P$ to $P(\mathbf{w})$. If $P(\mathbf{w})$ and $P(\mathbf{w}')$ are isomorphic, let $f:P(\mathbf{w})\to P(\mathbf{w}')$ be an isomorphism. Then $f$ is a linear extension of $P(\mathbf{w})$, and $\mathbf{w}'=\mathbf{w}(f)$, so $\mathbf{w}$ and $\mathbf{w}'$ are in the same commutation class. Conversely, if $\mathbf{w}$ and $\mathbf{w}'$ are in the same commutation class, let $e$ be a linear extension of $P(\mathbf{w})$ so that $\mathbf{w}'=\mathbf{w}(e)$. We then have that $e$ is an isomorphism between $P(\mathbf{w})$ and $P(\mathbf{w}')$.
\end{proof}

\section{Coxeter groups}

A group $W$ generated by a set $S$ (called the set of \emph{simple reflections}) is called a \emph{Coxeter group} if $W$ is determined by the relations $s^2=1$ for all $s\in S$, as well as zero or more relations of the form $(ss')^{m_{ss'}}=1$ for some integer $m_{ss'}\geq 2$, where $s,s'\in S$ and $s\neq s'$. We let $W$ be a Coxeter group and $S$ the chosen set of simple reflections from now on. A \emph{reduced word} for an element $w\in W$ is a word for $w$ of minimal length. We write $\ell(w)=\ell(\mathbf{w})$, where $\mathbf{w}$ is a reduced word for $w$.

For starters, we define a problem by the name of {\bf Reduced Word Count} as follows:\\

\emph{Input} A Coxeter group $W$ with generating set $S$ and an element $w\in W$.

\emph{Output} The number of reduced words for $w$.\\

Then we have
\begin{theorem}
The problem {\bf Reduced Word Count} is \#P-complete.
\end{theorem}
\begin{proof}
We can see that the problem is \#P-hard by noting that any poset occurs as a word poset for an element of a Coxeter group that has only one commutation class of reduced words. It follows from the numbers game (see \cite{BB}) that recognizing a reduced word for $w$ takes at most $O(n^2)$ time, where $n=l(w)$. Thus, {\bf Reduced Word Count} counts the number of accepting states of a Turing machine that runs in polynomial time, and hence is in \#P.
\end{proof}

We turn now to the study of commutation classes of reduced words in Coxeter groups. One benefit of restricting to this case is that there is a systematic way to find the commutation classes of reduced words of a given element.\\ 

We first need to define the \emph{root system} of a Coxeter group. These are \emph{not} necessarily the same as the root systems of Weyl groups from Lie theory, though they have many of the same properties. For a reference containing these facts, see, for example, \cite{BB}.

We recall that a Coxeter group is determined by the orders of the elements $ss'$ for $s,s'\in S$, which we denote by $m_{ss'}$. Thus $m_{ss'}$ is either a positive integer or $\infty$. We also have that $m_{ss'}=m_{s's}$, and $m_{ss}=1$ for all $s$. 

Let us assume that the generating set $S$ is finite. Let $V=\mathbb{R}^S$, the set of functions from $S$ to the real numbers. We choose a natural basis, namely $\{r_s\}_{s\in S}$, where $r_{s}(s)=1$ and $r_s(s')=0$ if $s'\neq s$. We define a bilinear form $(-,-)$ on $V$ by the rule
$$(r_{s},r_{s'})=\left\{\begin{array}{ll}
-1&\mbox{ if }m_{ss'}=\infty\\
-\cos{\frac{\pi}{m_{ss'}}}&\mbox{ otherwise.}\end{array}\right.$$
and extending linearly. Notice that this form is symmetric, and that $(r_{s},r_{s})=1$ for all $s$. There is a well-defined action of $W$ on $V$ given by the rule
$$sr_{s'}=r_{s'}-2(r_{s},r_{s'})r_{s}\mathrm{.}$$
With this action, $V$ is a representation of $W$, and the form $(-,-)$ is $W$-invariant: $(r,r')=(wr,wr')$ for all $w\in W$. 

A vector $r\in V$ is said to be a \emph{root} if $r=wr_{s}$ for some $s\in S$ and $w\in W$. The vectors $r_{s}$ are called the \emph{simple roots}, and the set of all roots is calld the \emph{root system} of $W$. The root system is a disjoint union of the sets $\Phi^+$ and $\Phi^-$, where $\Phi^+$ is the set of all roots $r$ such that the coefficient of $r_{s}$ is nonnegative for all $s$, and $\Phi^-$ is the set of all roots $r$ such that the coefficient of $r_{s}$ is nonpositive for all $s$. $\Phi^+$ is called the set of \emph{positive roots}, and $\Phi^-$ is the set of \emph{negative roots}.

For $w\in W$, a positive root $r$ is said to be an \emph{inversion} of $w$ if $wr$ is a negative root. We denote by $R(w)$ the set of inversions of $w^{-1}$. It is well-known that $|R(w)|=\ell(w)$; since $sr_{s}=-r_{s}$, this implies that if $r\in\Phi^+$ and $r\neq r_{s}$, then $sr\in\Phi^+$. We define the set $D_L(w)$ to be the set of all \emph{left descents} of $w$, meaning those $s\in S$ such that $l(sw)<\ell(w)$; then $r_{s}\in R(w)$ if and only if $s\in D_L(w)$. Analogously, $D_R(w)$, the set of \emph{right descents} of $w$, is the set of $s\in S$ such that $\ell(ws)<\ell(w)$. We will need the fact that if $s$ is not a right descent of $w$, then $R(w)\subset R(ws)$.\\

We will be studying functions $\lambda:R(w)\to\mathbb{N}_+$. If $\lambda:R(w)\to\mathbb{N}_+$ is a function and $s\notin D_R(w)$, then we define $\Phi(\lambda,s):R(ws)\to\mathbb{N}_+$ as follows. If $\lambda'=\Phi(\lambda,s)$, then $\lambda'(r)=\lambda(r)$ if $r\in R(w)$; if $r\notin R(w)$ and $r$ is simple, then $\lambda'(r)=1$; otherwise, $\lambda'(r)=\lambda(r')+1$, where $r'\in R(w)$ and $\lambda(r')$ is maximal with respect to the property that $(r,r')>0$.\\

For each $w\in W$, we define a set $C(w)$ of positive integer valued functions on the set $R(w)$ of inversions  of $w^{-1}$ as follows:\\

\begin{itemize}
\item[(1) ] If $e$ is the identity element, then $C(e)$ is the set containing the empty function.\\

\item[(2) ] If $\ell(w)>0$, then $C(w)$ contains exactly the functions $\Phi(\lambda,s)$ for $s\in D_R(w)$ and $\lambda\in C(ws)$.\\
\end{itemize}

%\begin{example}
%If $S=\{s_1,s_2\}$ and $m_{12}=3$, let $w=s_2s_1$. Then $R(w)=\{r_2,r_1+r_2\}$, and $C(w)$ contains the function $\lambda$ with $\lambda(r_2)=1$ and $\lambda(r_1+r_2)=2$. We have that $(r_1,r_2)=-\frac{1}{2}$, $(r_1,r_1+r_2)=\frac{1}{2}$, and $(r_2,r_1+r_2)=\frac{1}{2}$. To construct $\lambda'\in C(s_1w)$, we set $\lambda'(r_1)=1$, $\lambda'(s_1r_2)=\lambda'(r_1+r_2)=\lambda(r_2)+1=2$, and $\lambda'(s_1(r_1+r_2))=\lambda'(r_2)=\lambda(r_1+r_2)+1=3$.
%\end{example}

%The construction of the functions in $C(w)$ involves choosing a function $\lambda\in C(sw)$ for some $s\in D_L(w)$ and \emph{ascending} to obtain a unique corresponding $\lambda'\in C(w)$ with $\lambda'(r_{s})=1$. It is also possible to reverse this procedure and \emph{descend} from a $\lambda'\in C(w)$ with $\lambda'(r_{s})=1$ to a unique corresponding $\lambda\in C(sw)$ such that $\lambda'$ is obtained from $\lambda$ by ascending.\\

\newcommand{\dep}{\mathrm{dp}}

\noindent{}If $P$ is a finite partially ordered set, we define a function $\dep:P\to\mathbb{N}_+$ by declaring that $\dep(u)$ is the longest length of a strictly increasing sequence of elements $u_1<u_2<\cdots<u_{\dep(u)}=u$.\\

If $\mathbf{w}$ is a reduced word for an element $w\in W$, then there is a particular bijection $\phi_{\mathbf{w}}:P(\mathbf{w})\to R(w)$ that will be useful to us. For an index $i\in P(\mathbf{w})$, let us define
$$\phi_E(i)=\mathbf{w}_1\mathbf{w}_2\cdots\mathbf{w}_{i-1}r_{\mathbf{w}_i}\mathrm{.}$$

%\begin{example} \label{example:symroots}
%If $S=\{s_1,s_2\}$ and $m_{12}=3$, let $w=s_1s_2s_1$. Then $R(w)=\{r_1,r_1+r_2,r_2\}$. A reduced diagram $E$ for this element is
%$$\begin{array}{cc}
%1&2\\
%&1
%\end{array}$$
%Write this as
%$$\begin{array}{cc}
%x&y\\
%&z
%\end{array}$$
%Then $\phi_E(x)=r_1$, $\phi_E(y)=s_1r_2=r_1+r_2$, and $\phi_E(z)=s_1s_2r_1=r_2$. $\qed$
%\end{example}

\begin{lemma}
If $f:P(\mathbf{w})\to P(\mathbf{w}')$ is an isomorphism, then $\phi_{\mathbf{w}}=\phi_{\mathbf{w}'}\circ f$.
\end{lemma}
\begin{proof}
It suffices to show this when there is an $i$ such that $\mathbf{w}_j=\mathbf{w}'_j$ if $j\neq i,i+1$, $\mathbf{w}_i=\mathbf{w}'_{i+1}$, and $\mathbf{w}_{i+1}=\mathbf{w}'_i$. If $j<i$, then clearly $\phi_{\mathbf{w}}(j)=\phi_{\mathbf{w}'}(j)$. If $j>i+1$, then the element acting on $r_{\mathbf{w}_{j}}=r_{\mathbf{w}'_j}$ is the same, because the two words differ only by transposition of two adjacent commuting generators. $\mathbf{w}_ir_{\mathbf{w}_{i+1}}=r_{\mathbf{w}_{i+1}}$, so $\phi_{\mathbf{w}}(i+1)=\phi_{\mathbf{w}'}(i)$. Similarly, $\mathbf{w}_{i+1}r_{\mathbf{w}_i}=r_{\mathbf{w}_i}$, so $\phi_{\mathbf{w}}(i)=\phi_{\mathbf{w}'}(i+1)$.
\end{proof}

Thus, $\phi_{\mathbf{w}}$ depends only on the isomorphism class of $P(\mathbf{w})$. If $P$ is a word poset such that the corresponding commutation class consists of reduced words, we will call $P$ a \emph{reduced word poset}. Then $\phi_P$ is well-defined by the lemma.

\begin{lemma}
If $P$ is a reduced word poset for $w$, then $\phi_{P}:P\to R(w)$ is a bijection.
\end{lemma}
\begin{proof}
This is true if $\ell(w)=0$. If it holds for lengths less than $m$, suppose $\ell(w)=m$ and $\mathbf{w}$ is a reduced word for $w$. If $u$ is a minimal element of $P$, we have that $P'=P-\{u\}$ is a reduced word poset for the element $w'\in W$ with $w=s(u)w'$. We know that $R(w)=s(u)R(w')\cup \{r_{s(u)}\}$. Also, $\phi_{P}(s(u)r)=\phi_{P'}(r)$ for $r\in R(w')$, so $\phi_{P}|P'$ is a bijection by induction, and its image does not contain $r_{s(u)}$. Finally, $\phi_{P}(u)=r_{s(u)}$ by definition, so $\phi_{\mathbf{w}}$ is a bijection.
\end{proof}

The functions $\lambda\in C(w)$ will turn out to be the depth functions of the reduced word posets $P$ for $w$. For such a poset $P$, we define $\lambda_P:R(w)\to\mathbb{N}_+$ by $\lambda_P(r)=\dep(\phi_P^{-1}(r))$. 

%\begin{example}
%For the diagram $E$ in Example \ref{example:symroots}, we have $\lambda_E(r_1)=1$, $\lambda_E(r_1+r_2)=2$, and $\lambda_E(r_2)=3$. If we instead use the reduced diagram $E'$ given by
%$$\begin{array}{cc}
%2&\\
%1&2
%\end{array}$$
%and write this as
%$$\begin{array}{cc}
%a&\\
%b&c
%\end{array}$$
%then $\phi_{E'}(a)=r_2$, $\phi_{E'}(b)=s_2r_1=r_1+r_2$, and $\phi_{E'}(c)=s_2s_1r_2=r_1$, so that $\lambda_{E'}(r_1)=3$, $\lambda_{E'}(r_1+r_2)=2$, and $\lambda_{E'}(r_2)=1$. $\qed$
%\end{example}

Then we have

\begin{theorem} \label{theorem:depthfunctions}
$C(w)$ is identical to the set of functions $\lambda_P$, where $P$ ranges over all reduced word posets for $w$.
\end{theorem}

\begin{proof}
Let $P$ be a reduced word poset for $w$; we show that $\lambda_P\in C(w)$. If we choose a maximal element $u\in P$, then $u$ corresponds to a right descent of $w$. If we write $P'=P-\{u\}$, then by induction $\lambda_{P'}\in C(ws(u))$. If $\mathbf{w}'$ is a reduced word for $ws(u)$, then we may assume $P'=P(\mathbf{w}')$, and setting $\mathbf{w}=\mathbf{w}'s(u)$, we have $P(\mathbf{w})$. If $u$ was also a minimal element of $P$, then we are done. Otherwise, we may assume $s=\mathbf{w}'_{\ell(w)-1}$ and $s(u)$ do not commute, and let  $r'=xr_s$, and $r=xsr_{s(u)}$, where $s\in D_R(ws(u))$ and $x=wss(u)$. We may also assume that $\lambda_P(r')$ is maximal with respect to these equations. Then $(r',r)=(xr_s,xsr_{s'})=(r_s,sr_{s'})=-(r_s,r_{s'})>0$, so $\lambda_P(r)=\lambda_P(r')+1=\dep(u)$.

%By definition, $\lambda_P=\Phi(\lambda',s(u))$, so $\lambda_P\in C(w)$.

Suppose now that $\lambda\in C(w)$. We construct a diagram $P$ with $\lambda_P=\lambda$. Choose a root $r$ with $\lambda(r)$ maximal. Then $r$ corresponds to a right descent $s$ of $w$, so the function $\lambda'$ obtained by restricting to $R(w)-\{r\}=R(ws)$ is in $C(ws)$. By induction, there is a word poset $P'$ such that $\lambda'=\lambda_{P'}$. Then we may assume $P'=P(\mathbf{w}')$ for some reduced word $\mathbf{w}'$ for $ws$, and let us assume $\phi_{\mathbf{w}'}(\ell(\mathbf{w})-1)=r'$, where $\lambda_{P'}(r')=k-1$ and $(r,r')>0$. Write $\mathbf{w}=\mathbf{w}'s$, and let $P=P(\mathbf{w})$. Then $\lambda_P$ restricted to $P'$ is $\lambda_{P'}$. $\phi_{\mathbf{w}}(\ell(w))=r$, and $\phi_{\mathbf{w}}(\ell(w)-1)=r'$. We have that $r'=xr_s$, and $r=xsr_{s'}$. Then $(r,r')=(xr_s,xsr_{s'})=(r_s,sr_{s'})=-(r_s,r_{s'})>0$. Thus, $s$ and $s'$ do not commute, so $\dep(\phi_{w}^{-1}(r))=\dep(\phi_{w}^{-1}(r'))+1=k$.
\end{proof}

As a corollary, we have
\begin{corollary}
$C(w)$ is in bijective correspondence with the commutation classes of reduced words of $w$.
\end{corollary}
\begin{proof}
Let $P,P'$ be reduced word posets for $w$. It is clear from the definition of $\lambda_P$ and the lemma that if $P$ and $P'$ are isomorphic, then $\lambda_P=\lambda_{P'}$. We can also see by induction that if $P$ and $P'$ are not isomorphic, then $\lambda_P\neq \lambda_{P'}$.
\end{proof}

%We present now an extended example. Suppose $S=\{s_1,s_2,s_3\}$, $m_{12}=m_{23}=3$, and $s_1s_3=s_3s_1$ (i.e. $m_{13}=2$). $W$ is then isomorphic to $S_4$, the symmetric group on 4 letters, and the longest element of $W$ is  $w_0=s_1s_2s_1s_3s_2s_1$.  Below is the root system of $W$, where two roots $r,r'$ are connected by a line in the diagram if $(r,r')\neq 0$.

%{\includegraphics{roots.eps}}

%The 8 functions in $C(w_0)$, with corresponding diagrams, are then

%{\includegraphics{c1.eps}}
%\begin{tabular}{ccc}
%1&2&3\\
%&1&2\\
%&&1
%\end{tabular}

%
%\end{pspicture} 

%{\includegraphics{c2.eps}}
%\begin{tabular}{ccccc}
%&&3&&\\
%&&2&&\\
%&&1&2&3\\
%&&&  &2
%\end{tabular}

%{\includegraphics{c3.eps}}
%\begin{tabular}{ccccc}
%&&3&&\\
%&&2&3&\\
%&&1&2&3
%\end{tabular}

%132543

%{\includegraphics{c4.eps}}
%\begin{tabular}{ccccc}
%&2&3&&\\
%& &2&&\\
%&& 1&2&3
%\end{tabular}

%{\includegraphics{c5.eps}}
%\begin{tabular}{ccccc}
%&2&&&\\
%&1&2&3&\\
%& & &2&\\
%& & &1&
%\end{tabular}

%{\includegraphics{c6.eps}}
%\begin{tabular}{ccccc}
%&2&3&&\\
%&1&2&3&\\
%& &1&&
%\end{tabular}

%{\includegraphics{c7.eps}}
%\begin{tabular}{ccccc}
%& &3&&\\
%&1&2&3&\\
%& &1&2&
%\end{tabular}

%{\includegraphics{c8.eps}}
%\begin{tabular}{ccccc}
%1&2&3&&\\
%&  &2&&\\
%& & 1&2&
%\end{tabular}

It now becomes of interest to compute the size of $C(w)$. We present here a recurrence relation. We say a subset $T\subset S$ is \emph{independent} if all elements of $T$ commute with each other. If $T$ is an independent subset of $S$, we may consider $T$ as an element of $W$ by defining $T=\prod_{s\in T}{s}$. This is well-defined because all elements of $T$ commute. First, a lemma.

\begin{lemma}
For all $\lambda_P\in C(w)$ and all $k\in\mathbb{N}_+$, the elements in $P_k=\phi_P^{-1}(\lambda_P^{-1}(\{k\}))$ are pairwise incomparable, and hence $s(P_k)$ is an independent subset of $G$ such that $|s(P_k)|=|P_k|$.
\end{lemma}
\begin{proof}
$\phi_P^{-1}(\lambda_P^{-1}(\{k\}))$ is a set of elements of $P$, all of which have the same depth, so the elements are pairwise incomparable. $|s(P_k)|=|P_k|$ because two elements with the same label would be comparable, and $s(E_k)$ is independent because two elements whose labels did not commute would be comparable.
\end{proof}

\begin{theorem}
Let $w\in W$. Then we have
$$|C(w)|=\sum_{\substack{\emptyset\neq T\subset D_R(w)\\T\mbox{ independent}}}{(-1)^{|T|+1}|C(wT)|}\mathrm{.}$$
\end{theorem}

\begin{proof}
For $T\subset D_R(w)$, write $C_A(w,T)$ to mean the set of all $\lambda\in C(w)$ such that $\lambda(r_{s_i})=1$ for all $s_i\in T$. We also write $C_=(w,T)$ for $T'\subset D_R(w)$ to mean the set of all $\lambda\in C(w)$ such that $\lambda(r_{s_i})=1$ if and only if $s_i\in T'$. Note that $C_A(w,T)$ and $C_=(w,T)$ are empty if $T$ is not independent, by the lemma. We have
$$C_A(w,T)=\sum_{T\subset T'}{C_=(w,T')}$$
so, by inclusion-exclusion,
$$C_=(w,T)=\sum_{T\subset T'}{(-1)^{|T'-T|}C_A(w,T')}$$
and summing over all independent subsets of $D$ we get
$$|C(w)|=\sum_{T\subset D}{\sum_{T\subset T'}{(-1)^{|T'-T|}C_A(w,T')}}\mathrm{.}$$
Interchanging the order of summation, we get
$$|C(w)|=\sum_{T'\subset D}{\left(\sum_{T\subset T'}{(-1)^{|T'-T|}}\right)C_A(w,T')}$$
and that works out to
$$|C(w)|=\sum_{T'\subset D}{(-1)^{|T'|+1}C_A(w,T')}$$
But if $T'\subset D$ is independent and nonempty, $C_A(w,T')\subset C(w)$ can be naturally identified with the contribution from $C(T'w)$.
\end{proof}

Since $|C(w)|=|C(w^{-1})|$, we have

\begin{corollary}
Let $w\in W$ and write $D=D_R(w)=D_L(w^{-1})$. Then
$$|C(w)|=\sum_{\substack{\emptyset\neq T\subset D\\T\mbox{ independent}}}{(-1)^{|T|+1}|C(wT)|}\mathrm{.}\qed$$
\end{corollary}

\begin{example}
Let us use this formula to compute $|C([4231])|$. A right descent of a permutation is the same as a pair of adjacent elements such that the first is greater than the second, and a set of descents is independent if and only if the first elements of the pairs are each at least two spaces apart from each other. We write $C[abcd]$ to mean $|C([abcd])|$ for brevity's sake, and to shorten the computation we use the fact that $C(w)=1$ if there are no indices $i<j<k$ such that $w(i)>w(j)>w(k)$. Then
$$C[4231]=C[2431]+C[4213]-C[2413]=C[2431]+C[4213]-1$$
$$C[2431]=C[2341]+C[2413]=1+1=2$$
$$C[4213]=C[2413]+C[4123]=1+1=2$$
so
$$C[4231]=2+2-1=3\mathrm{.}$$
\end{example}

\end{document}